\documentclass[10pt]{amsart}
\usepackage{verbatim, graphics, color}
\usepackage{tikz,pgf}

\newtheorem{theorem}{Theorem}
\newtheorem{proposition}[theorem]{Proposition}

\newtheorem{corollary}[theorem]{Corollary}
\newtheorem{conjecture}[theorem]{Conjecture}

\theoremstyle{definition}
\newtheorem{definition}{Definition}

\newcommand{\commsq}[8]{ \node (a) at (-1,1) {#1};
 \node (b) at (1,1) {#2};
 \node (c) at (1,-1) {#3};
 \node (d) at (-1,-1) {#4};
 \draw[->] (a) -- node[above] {#5} (b);
 \draw[->] (b) -- node[right] {#6} (c);
 \draw[->] (a) -- node[left] {#8} (d);
 \draw[->] (d) -- node[above] {#7} (c);
}

\newcommand{\spem}[1]{ \mathbf{#1}}
\newcommand{\spemg}[2]{ \spem{#1}[#2]}
\newcommand{\spemi}[3]{\spemg{#1}{#2_{#3}}}
\newcommand{\cat}[1]{ \mathsf{#1}}
\newcommand{\spe}{\cat{Sp}}
\newcommand{\gs}{\cat{GrSp}}
\newcommand{\id}{\mbox{ id}}

\theoremstyle{remark}
\newtheorem{examples}{Examples}

\newtheorem{remark}[theorem]{Remark}
\newtheorem{example}[theorem]{Example}

\author{Jacob A. White }
\title{On some Hopf monoids in Graphical Species}
\address{ School of Mathematical and Statistical Sciences, Arizona State University, Tempe, AZ}
\keywords{Species, Combinatorial Hopf Algebras, Graph Theory, Hyperplane Arrangements}

\begin{document}

\maketitle

\begin{abstract}

Combinatorial Hopf algebras arise in a variety of applications. Recently, Aguiar and Mahajan showed how many well-studied Hopf algebras are closely related to Hopf monoids in species. 
 In this paper, we study Hopf monoids in graphical species, giving a `graph-theoretic' analogue to the work of Aguiar and Mahajan. In particular, several examples of Hopf monoids in graphical species are detailed, most of which are related to graph coloring, or hyperplane arrangements associated to graphs.

\end{abstract}

\section{Introduction}

The notion of combinatorial bialgebra goes back to Joni and Rota \cite{joni-rota}. The idea was to use terminology from algebra to explain the relationship between a way of combining and decomposing combinatorial objects. The idea was that combining two combinatorial objects gave a product, and decomposing gave a coproduct, and the two were related by a braiding axiom.

The notion of combinatorial Hopf algebra has been extended to Joyal's category of species by Aguiar and Mahajan \cite{thebook}, so we can now speak of Hopf monoids in species. The idea is that now our combinatorial objects come with labels (from some finite set $S$), and the objects do not depend on the labels. In this case, the notion of product is to combine combinatorial objects whose label sets are disjoint, and the notion of coproduct is to decompose a combinatorial object into pairs of objects with disjoint label sets.

One can discuss studying other monoidal categories of interest in combinatorics. Bergeron and Choquette \cite{hyperoctahedral} have studied $\mathcal{H}$-species, which come from hyperoctahedral groups. They also studied various functors, and the resulting Hopf algebras.
In some sense, they were looking at combinatorial objects of `type $B$'. For the present paper, we are instead interested in graph-theoretic objects.

In some sense, the work of Aguiar and Mahajan \cite{thebook} gives rise to the idea that, instead of attempting to view something as being a combinatorial Hopf algebra, it makes sense to consider what category the combinatorial objects really belong to. For instance, labeled combinatorial objects should be viewed as species. In many cases, labeled combinatorial objects can be generalized to being combinatorial objects on graphs (and usually one recovers the original labeled combinatorial objects by restricting to complete graphs). For instance, linear orders generalize to acyclic orientations, labeled trees generalize to spanning trees, and so on.  
Hence, there is some motivation to asking whether or not the corresponding `graph-theoretic' analogues also form Hopf monoids.

In this paper, we initialize the study of Hopf monoids in the category of graphical species. That is, now we are studying combinatorial structures on finite graphs, and various ways to combine and decompose them across induced subgraphs. It turns out that many well-known combinatorial structures (matchings, stable partitions, acyclic orientations) carry Hopf monoid structures when viewed as graphical species. Moreover, some of these graphical species come with Hopf monoid morphisms, and the Hopf monoid structures are a `graph-theoretic' analogue of known Hopf monoids in species. 

This paper  will mainly serve for giving definitions of graphical species, as well as detailing several examples. In future papers, we shall discuss other interesting ideas motivated by this paper, and related to graphical species.

The layout is as follows: in the next section we review the definition of the category $\spe$ of species, its lax braided monoidal structure, and several examples of Hopf monoids in species. We also review the notion of bilax monoidal functor. In Section 3, we define the category of graphical species, $\gs$, give it a lax braided monoidal category, and define Hopf monoids in $\gs$ with respect to the braided monoidal structure. In Section 4, we study examples of Hopf monoids related to linear orders, including acyclic orientations and stable compositions. We also relate these Hopf monoids via commutative diagrams. In Section 5, we study examples related to set partitions, including stable partitions and flats. We also relate these Hopf monoids through commutative diagrams. In Section 6, we study how graph complementation gives rise to a bistrong endofunctor on $\gs$ that changes the braiding map. Then we mention some bilax monoidal functors from $\gs$ to $\spe$, which explain why several of the known Hopf monoids in $\spe$ are `shadows' of the Hopf monoids in $\gs$ we discuss in this paper. Finally, in Section 7, we mention several future directions.



\section{Preliminaries}

\subsection{Species}

Throughout this paper, fix a field $\mathbb{K}$.
For a definition of lax braided monoidal category, see Joyal and Street \cite{joyal-street} or Aguiar and Mahajan \cite{thebook}. For the definitions of monoids, comonoids, bimonoids and Hopf monoids in monoidal categories, see Aguiar and Mahajan \cite{thebook}. We will discuss these definitions in the context of the lax braided monoidal category of graphical species with respect to Cauchy product.

Let $\cat{Set}$ be the category of all finite sets with bijections as morphisms.
Let $\cat{Vec}$ be the category of all vector spaces over $\mathbb{K}$ with linear maps as morphisms.

 A \emph{species} is a functor $\spem{p}: \cat{Set} \to \cat{Vec}$. That is, for each finite set $I$, we associate a vector space $\spem{p}[I]$, and to each bijection $\sigma: I \to J$, we associate a linear map $\spem{p}[\sigma]: \spem{p}[I] \to \spem{p}[J]$, such that 
$\spem{p}[id_I] = id_{\spem{p}[I]}$, and $\spem{p}[\sigma \circ \tau] = \spem{p}[\sigma] \circ \spem{p}[\tau]$ whenever $\sigma \circ \tau$ is defined.

 A morphism $\varphi: \spem{p} \to \spem{q}$ of species is a natural transformation from $\spem{p}$ to $\spem{q}$. That is, for each finite set $I$ we have a linear map $\varphi_I: \spem{p}[I] \to \spem{q}[I]$ such that for any bijection $\sigma: I \to J$ we have the following commutative diagram:
 \begin{figure} [h]
\label{diag:speciesmorphism}
\begin{center}
\begin{tikzpicture}
 \commsq{$\spem{p}[I]$}{$\spem{q}[I]$}{$\spem{q}[J]$}{$\spem{p}[J]$}
{$\varphi_I$}{$\spem{q}[\sigma]$}{$\varphi_J$}{$\spem{p}[\sigma]$}

 \end{tikzpicture}
\end{center}
\end{figure}

Let $\spe$ denote the category of species with natural transformations as morphisms. The category of species was first introduced by Joyal \cite{joyal}. The book by Bergeron, Labelle and Leroux \cite{species-book} forms a wonderful introduction into combinatorial species.

\begin{examples} Here are some examples of species.
\begin{enumerate}

\item Let $\spem{E}$ be the species defined by $\spem{E}[I] = \mathbb{K}$ for every finite set $I$. For an isomorphism $\sigma: I \to J$, let $\spem{E}[\sigma] = id_{\spem{E}[I]}$. $\spem{E}$ is called the \emph{exponential} species.

\item Let $\spem{L}$ be given by defining $\spem{L}[I]$ to be the vector space over $\mathbb{K}$ with basis given by all linear orders on $I$.

\item Let $\spem{\Pi}$ be given by defining $\spem{\Pi}[I]$ to be the vector space over $\mathbb{K}$ with basis given by all set partitions of $I$.

\item Let $\spem{\Sigma}$ be given by defining $\spem{\Pi}[I]$ to be the vector space over $\mathbb{K}$ with basis given by all set compositions (ordered set partitions) of $I$.
\end{enumerate}
\end{examples}

We review the Cauchy product, and substitution product for species.
 Given a finite set $S$, a composition of $S$, $C = S_1 | S_2 | \ldots | S_k$, is a sequence of disjoint subsets of $S$ whose union is $S$. We often write $C \models S$ to say that $C$ is a composition of $S$.
A set partition $\pi = \{S_1, \ldots, S_k \}$ of $S$, is a collection of disjoint sets whose union is $S$. We often write $\pi \vdash S$ to say that $\pi$ is a partition of $S$. Also, we refer to the elements of $\pi$ as blocks.

 Let $\spem{p}$ and $\spem{q}$ be species.
We define the $I$-component of their Cauchy product by:
\begin{equation}
 \label{eq:speciescauchy}
(\spem{p} \cdot \spem{q})[I] = \bigoplus_{S|T \models I} \spem{p}[S] \otimes \spem{q}[T]
\end{equation}

If $q$ is a positive species, then we define the $I$-component of the substitution product by:
\begin{equation}
 \label{eq:speciescomposition}
(\spem{p} \circ \spem{q})[I] = \bigoplus_{\pi \vdash I} \spem{p}[\pi] \otimes \bigotimes_{B \in \pi} \spem{q}[B]
\end{equation}

First, we recall how to turn $\spe$ into a lax braided monoidal category \cite{thebook}. Fix $q \in \mathbb{K}$.
We recall the braiding map $\beta_q: \spem{p} \cdot \spem{q} \rightarrow \spem{q} \cdot \spem{p}$. The $I$-component of the braiding, is the direct sum over all compositions $S|T \models I$, of the maps
$$\beta_{q,S,T}: \spem{p}[S] \otimes \spem{q}[T] \rightarrow \spem{q}[T] \otimes \spem{p}[S]$$
given by $\beta_{q,S,T}(x \otimes y) = q^{|S||T|} y \otimes x$.
It is known that $(\spe, \cdot, \beta_q)$ forms a lax braided monoidal category.

In particular, there is a notion of Hopf monoid in species. In Chapter 12 of their book, Aguiar and Mahajan \cite{thebook} construct several Hopf monoids in species related to the braid arrangement. These include Hopf monoids on $\spem{L}$, $\spem{\Pi}$, $\spem{E}$, and $\spem{\Sigma}$.

The main goal of this paper is to extend this work, by moving from species to graphical species. In particular, we demonstrate several graphical species which generalize the above species, and study Hopf monoid structures on these graphical species. We also investigate several bilax and bistrong monoidal functors coming from graphical species.

\subsection{Bilax monoidal functors}
\begin{definition}[B\'enadou \cite{benabou} and Aguiar and Mahajan \cite{thebook}]
Let $(C,\bullet,e)$ and $(D,\diamond,e')$ be two monoidal categories.
\begin{enumerate}
\item A lax monoidal functor $(F,\phi,\phi_0)$ between $(C,\bullet,e)$ and $(D,\star,e')$ is a functor $F\colon C\rightarrow D$, and a morphism $\phi_{a,b}\colon F(a) \star F(b) \rightarrow F(a\bullet b)$ in $D$, which is natural in $a$ and $b$ for each pair $a,b$ of objects in $C$, and a morphism $\phi_0\colon e' \rightarrow F(e)$ in $D$ such that $\phi$ is associative and left and right unital in the usual sense.
\item A lax monoidal functor $(F, \phi, \phi_0)$ between lax braided monoidal categories $(C, \bullet, e, \beta)$ and $(D, \star, e', \beta')$ is lax braided if $\phi_{b,a} \circ \beta'_{F(a),F(b)}) = F(\beta_{a, b}) \circ \phi_{a,b}$.
\item A colax monoidal functor $(F,\psi,\psi_0)$ is a functor $F\colon C\rightarrow D$ with a morphism $\phi_{a,b}\colon F(a\bullet b) \rightarrow F(a) \star F(b)$ in $D$, which is natural in $a$ and $b$ for each pair of objects $a$, $b$ in $C$, and a morphism $\psi_0\colon F(e) \rightarrow e'$ in $D$ such that $\psi$ is coassociative and left and right counital in the usual sense.
\item A colax monoidal functor $(F, \psi, \psi_0)$ between lax braided monoidal categories $(C, \bullet, e, \beta)$ and $(D, \star, e', \beta')$ is colax braided if $\psi_{b,a} \circ F(\beta_{a,b}) = \beta'_{F(a), F(b)} \circ \psi_{a,b}$.
\item A bilax monoidal functor $(F,\phi,\phi_0,\psi,\psi_0)$ between two braided monoidal categories $(C,\bullet,e,\beta)$ and $(D,\diamond,e',\beta)$ is a lax monoidal functor $(F,\phi,\phi_0)$ and a colax monoidal functor $(F,\psi,\psi_0)$ satisfying the braiding condition and unitality conditions.
\item A bilax monoidal functor with $\phi$, $\phi_0$, $\psi$ and $\psi_0$ invertible is a bistrong monoidal functor.
\end{enumerate}
\end{definition}

The composites of lax, colax, bilax and bistrong monoidal functors are lax, colax, bilax and bistrong monoidal functors respectively.
A morphism of bilax monoidal functors between $(F,\phi,\psi)$ and $(G,\gamma,\delta)$ is a natural transformation $\alpha \colon F\rightarrow G$ which commutes with $\phi$, $\psi$, $\gamma$ and $\delta$. See \cite{thebook} for more details.

\begin{proposition}[B\'enabou \cite{benabou} and Aguiar and Mahajan \cite{thebook}]\label{prop:Fhbimonoid}
\begin{enumerate}

\item If $F$ is a lax (colax, bilax) monoidal functor from $(C,\bullet)$ to $(D,\diamond)$ and $h$ is a monoid (co\-mon\-oid, bimonoid) in $C$ then $F(h)$ is a monoid (co\-mon\-oid, bimonoid) in $D$. Furthermore, if $F$ is braided lax (braided colax), and $h$ is commutative (cocommutative), then $F(h)$ is also commutative (cocommutative).
\item If $F$ is a bistrong monoidal functor from $(C,\bullet)$ to $(D,\diamond)$ and $h$ is a Hopf monoid in $C$ with antipode $S$ then $F(h)$ is a Hopf monoid in $D$ with antipode $F(S)$.
\end{enumerate}
\end{proposition}

\section{The notion of Hopf monoid for graphical species}

\subsection{Graphical species}
Let $\cat{Graph}$ be the category of all finite graphs, with graph isomorphisms as morphisms.
 A \emph{graphical species} is a functor $\spem{g}: \cat{Graph} \to \cat{Vec}$. That is, for each finite graph $G$, we associate a vector space $\spem{g}[G]$, and to each graph isomorphism $\sigma: G \to H$, we associate a linear map $\spem{g}[\sigma]: \spem{g}[G] \to \spem{g}[H]$, such that 
$\spem{g}[id_G] = id_{\spem{g}[G]}$, and $\spem{g}[\sigma \circ \tau] = \spem{g}[\sigma] \circ \spem{g}[\tau]$ whenever $\sigma \circ \tau$ is defined.

In other words, a graphical species is a vector-space valued graph invariant. Numerical and polynomial graph invariants are often studied in the literature. It is not hard to see that most graph invariants could be approached using graphical species instead. 

 A morphism $\varphi: \spem{g} \to \spem{h}$ of graphical species is a natural transformation from $\spem{g}$ to $\spem{h}$. That is, for each finite graph $G$ we have a linear map $\varphi_G: \spem{g}[G] \to \spem{h}[G]$ such that for any graph isomorphism $\sigma: G \to H$ we have the following commutative diagram:
\begin{center}
\begin{tikzpicture}
 \commsq{$\spem{g}[G]$}{$\spem{h}[G]$}{$\spem{h}[H]$}{$\spem{g}[H]$}
{$\varphi_G$}{$\spem{h}[\sigma]$}{$\varphi_H$}{$\spem{g}[\sigma]$}

 \end{tikzpicture}
\end{center}

Let $\gs$ denote the category of graphical species with natural transformations as morphisms.

\begin{example}
 We define the exponential graphical species $\spem{E}$ as follows: $\spem{E}[G] = \mathbb{K}$ for every finite graph $G$. We define 
$\spem{E}[\sigma] = id_{\spem{E}[G]}$ for any graph isomorphism $\sigma: G \to H$.
\end{example}

We say a (graphical) species is connected if $\spem{p}[\emptyset] = \mathbb{K}$, and we say the (graphical) species is positive if $\spem{p}[\emptyset] = 0$.

\subsection{Cauchy product and braiding for graphical species}

Now we define the Cauchy product for graphical species. We invite the reader to note the similarity to the Cauchy product for species.
Let $\spem{g}$ and $\spem{h}$ be graphical species.
Given a graph $G$, and a subset $S \subset V(G)$, let $G_S$ denote the induced subgraph (which has vertex set $S$ and edge set $\{uv: uv \in E(G), u, v \in S \}$). 
Also, recall the notion of \emph{quotient graph}. Given a vertex partition $\pi \vdash V(G)$, the graph $G / \pi$ has vertex set $\{B: B \in \pi \}$, and edges $BC$ if and only if there exists $b \in B, c \in C$ such that $bc \in E(G)$. 
We define the $G$-component of the Cauchy product by
\begin{equation}
 \label{eq:product}
(\spem{g} \cdot \spem{h})[G] = \bigoplus_{S|T \models V(G)} \spem{g}[G_S] \otimes \spem{h}[G_T]
\end{equation}

If $\spem{h}$ is a positive graphical species, we define the $G$-component of the substitution product by:
\begin{equation}
 \label{eq:composition}
(\spem{g} \circ \spem{h})[G] = \bigoplus_{\pi \vdash I} \spem{g}[G/\pi] \otimes \bigotimes_{B \in \pi} \spem{h}[G(B)]
\end{equation}

Now we turn $\gs$ into a lax braided monoidal category. Fix $q,t \in \mathbb{K}$. Given a graph $G$, let $\overline{G}$, the complementary graph, have vertex set $V(G)$ and edges $uv$ if and only if $uv \not\in E(G)$.
Given two graphical species, $\spem{g}$ and $\spem{h}$, a finite graph $G$, and a composition $S|T \models V(G)$, define $\beta_{q,t,G}^{S,T}: \spemi{g}{G}{S} \otimes \spemi{h}{G}{T} \to \spemi{h}{G}{T} \otimes \spemi{g}{G}{S}$ by $\beta_{q,t,G}^{S,T}(x \otimes y) = q^{e(G,S,T)} t^{e(\overline{G},S,T)} y \otimes x$, where $e(G,S,T)$ is the number of edges of $G$ with one endpoint in $S$ and one endpoint in $T$.  
Let $\beta_{q,t,G}$ be given by 
\begin{equation}
 \label{eq:braiding}
\beta_{q,t,G} = \oplus_{S|T \models V(G)} \beta_{q,t,G}^{S,T}
\end{equation}
Some particular choices of $q,t$ will be nice to refer to. We let $\beta = \beta_{1,1}$, $\beta_q = \beta_{q,1}$, and $\overline{\beta}_q = \beta_{1,q}$.

\begin{proposition}
 $(\gs, \cdot, \beta_{q,t})$ is a lax braided monoidal category, with product given by (\ref{eq:product}), with braiding 
$\beta_{q,t}$ whose $G$-component is given by (\ref{eq:braiding}), and unit given by 
$$\spem{o}[G] = \left\{ \begin{array}{ll} \mathbb{K} & \mbox{if } G = \emptyset \\ 0 & \mbox{otherwise} \end{array} \right.$$
\end{proposition}
The axioms are straightforward to check. Note that when $q,t \in \{+1, - 1 \}$, $(\gs, \cdot, \beta_{q,t})$ is a symmetric monoidal category, and whenever $q,t\neq 0$ $(\gs, \cdot, \beta_{q,t})$ is braided. 


\subsection{$(q,t)$-Hopf monoids in graphical species}

We now proceed to define the terms monoid, comonoid, bimonoid, and Hopf monoid for graphical species, with respect to the Cauchy product.
A \emph{monoid} in $\gs$ is a graphical species $\spem{g}$ together with a multiplication map $\mu$ and unit map $\eta$:
$$\mu: \spem{g} \cdot \spem{g} \to \spem{g}, \eta:\spem{o} \to \spem{g}$$
which are associative and unital in the usual sense. In particular, for each finite graph $G$ and composition $S|T \models V(G)$, we have a linear map
$$\mu_G^{S,T}: \spemi{p}{G}{S} \otimes \spemi{p}{G}{T} \to \spemg{p}{G}$$
and one addition linear map for the unit
$$\eta_{\emptyset}: \mathbb{K} \to \spemg{p}{\emptyset}.$$
Note that the $G$-component of $\mu$ is given by $\mu_G = \oplus_{S|T \models V(G)} \mu_G^{S,T}$.

A \emph{comonoid} in $\gs$ is a graphical species $\spem{g}$ together with a comulitplication map $\Delta$ and unit map $\epsilon$:
$$\Delta: \spem{g} \to \spem{g} \cdot \spem{g}, \epsilon:\spem{g} \to \spem{o}$$
which are coassociative and counital in the usual sense. In particular, for each finite graph $G$ there is a linear map $$\Delta_G: \spem{g}[G] \to \oplus_{S|T \models V(G)} \spemi{p}{G}{S} \otimes \spemi{p}{G}{T}$$
and for each each composition $S|T \models V(G)$, we obtain a linear map
$$\Delta_G^{S,T}: \spem{g}[G] \to \spemi{p}{G}{S} \otimes \spemi{p}{G}{T}.$$
Finally, we have one linear map for the counit:
$$\epsilon_{\emptyset}: \spemg{p}{\emptyset} \to \mathbb{K}.$$

A $(q,t)$-\emph{bimonoid} is a species $\spem{g}$ with monoid and comonoid structures such that the monoid and comonoid structures are compatible via the usual braiding axioms, where the braiding map is the map $\beta_{q,t}$. We refer to bimonoids in $(\gs, \cdot, \beta_q)$ as $q$-bimonoids, and bimonoids in $(\gs, \cdot, \beta)$ as bimonoids. In the case $q,t \neq 0$, one can equivalently define a $(q,t)$-bimonoid to be a graphical species $\spem{g}$, with monoid and comonoid structure, such that $\Delta$ and $\epsilon$ are monoid morphisms. 
We do not study bimonoids in $(\gs, \cdot, \bar{\beta}_q)$ in this paper.

A $(q,t)$-\emph{Hopf monoid} is a $(q,t)$-bimonoid $\spem{g}$ with a morphism $s:\spem{g} \to \spem{g}$, the \emph{antipode}, such that for each graph non-empty graph $G$ and vertex decomposition $A|B \models V(G)$, we have that the following composite maps
\begin{equation}
\label{eq:antipode}
\begin{tikzpicture}[>= latex]
\node (a) at (0,0) {$\spem{g}[G]$};
\node (b) at (2.5,0) {$\spem{g}[G_A] \otimes \spem{g}[G_B]$};
\node (c) at (6.5,0) {$\spem{g}[G_A] \otimes \spem{g}[G_B]$};
\node (d) at (9,0) {$\spem{g}[G]$};
\draw[->] (a) -- node[above] {$\Delta_G^{A,B}$} (b);
\draw[->] (b) -- node[above] {$s_{G_A} \otimes \id_{G_B}$} (c);
\draw[->] (c) -- node[above] {$\mu_G^{A,B}$} (d);
\end{tikzpicture}
\end{equation}
\begin{displaymath}
\begin{tikzpicture}[>= latex]
\node (a) at (0,0) {$\spem{g}[G]$};
\node (b) at (2.5,0) {$\spem{g}[G_A] \otimes \spem{g}[G_B]$};
\node (c) at (6.5,0) {$\spem{g}[G_A] \otimes \spem{g}[G_B]$};
\node (d) at (9,0) {$\spem{g}[G]$};
\draw[->] (a) -- node[above] {$\Delta_G^{A,B}$} (b);
\draw[->] (b) -- node[above] {$\id_{G_A} \otimes s_{G_B}$} (c);
\draw[->] (c) -- node[above] {$\mu_G^{A,B}$} (d);
\end{tikzpicture}
\end{displaymath}
are zero, and
\begin{displaymath}  
\mu_{\emptyset} \circ (\id_{\emptyset} \otimes s_{\emptyset}) \circ \Delta_{\emptyset} = \eta_{\emptyset} \circ \epsilon_{\emptyset} = \mu_{\emptyset} \circ (s_{\emptyset} \otimes \id_{\emptyset}) \circ \Delta_{\emptyset}.
\end{displaymath}

We say a monoid $(\spem{g}, \mu, \epsilon)$ is \emph{commutative} if $\mu = \mu \circ \beta_{q,t}$, and a comonoid $(\spem{g}, \Delta, \eta)$ is cocommutative if $\Delta = \beta_{q,t} \circ \Delta$. For graphical species, we have two new algebraic properties: disjoint commutativity, and join commutativity. We say that a monoid $(\spem{g}, \mu, \epsilon)$ is \emph{disjoint commutative} if $\mu_G^{S,T} = \mu_G^{T,S} \circ \beta_{q,t, G}^{S,T}$ whenever $e(G,S,T) = 0$. We say that a monoid is \emph{join commutative} if $\mu_G^{S,T} = \mu_G^{T,S} \circ \beta_{q,t,G}^{S,T}$ whenever $e(\overline{G}, S,T) = 0$. That is, for graphical species we have monoids that are commutative when the product is taken across disconnected subgraphs, or totally connected subgraphs.

\subsection{Takeuchi and Milnor-Moore formulas}
Much like in the case of species, if $\spem{g}$ is a connected $(q,t)$-bimonoid, then it is a $(q,t)$-Hopf monoid, and there is a formula for the antipode $s$. 
The proof is similar to Aguiar and Mahajan \cite{thebook}.

To state the formula, recall that since we are working in a braided monoidal category, given a monoid $\spem{g}$, and a positive integer $k > 2$, there is always a map $\mu_k: \spem{g}^{\cdot k} \to \spem{g}$ which is the $k$-fold multiplication. Let $\mu_G^{S_1, \ldots, S_k}$ denote the component of the $k$-fold multiplication coming from the finite graph $G$ with composition $S_1|\ldots|S_k \models V(G)$. Likewise, for a comonoid, there is a natural $k$-fold comultiplication $\Delta_k: \spem{P}^{\cdot k} \to \spem{g}$, and we let $\Delta_G^{S_1, \dots, S_k}$ denote the component of the comultiplication coming from the finite graph $G$ with vertex composition $S_1|\ldots|S_k \models V(G)$.
\begin{theorem}[Takeuchi's Formula]
 Let $\spem{g}$ be a connected $(q,t)$-bimonoid. Then $\spem{g}$ is a $(q,t)$-Hopf monoid, with antipode defined by:
$$s_{\emptyset} = 0$$
If $G$ is non-empty, then $$s_G = \sum_{k \geq 0} (-1)^k \sum_{S_1|\ldots|S_k \models V(G)} \mu_G^{S_1, \ldots, S_k} \circ \Delta_G^{S_1, \ldots, S_k}$$
where we are summing over all compositions $S_1|\ldots|S_k$ where the $S_i$ are non-empty.
\end{theorem}
Note that the above formula, for a typical bimonoid, may have many cancellations. We often study Hopf monoids where the graphical species $\spem{g}$ is equipped with a basis indexed by combinatorial structures on graphs. In these cases, we find simpler, combinatorial expressions for the antipode. 

For our purposes, there is also a recursive formula for the antipode, which is helpful for doing inclusion-exclusion. 
\begin{theorem}[Milnor-Moore Formula]
 Let $\spem{g}$ be a connected $(q,t)$-Hopf monoid. Then the antipode $s$ is defined by:
$s_{\emptyset} = 0$
If $G$ is non-empty, then $$s_G = - \sum_{S|T \models V(G), S \neq \emptyset} \id_{G_S} \otimes s_{G_T} $$ or 
$$s_G = - \sum_{S|T \models V(G), T \neq \emptyset} \id_{G_S} \otimes s_{G_T}$$
\end{theorem}

\section{Graphical Species generalizing $\spem{L}$}
In this section, we construct three graphical species that generalize the species of linear orders.
The first one is $\spem{L}_{q,t}$, the graphical species of linear orders. The second one is $\overrightarrow{\spem{A}}_q$, the graphical species of acyclic orientations. The final one is 
$\spem{S\Sigma}_{q,t}$, the graphical species of stable compositions.
 
\subsection{The graphical species of linear orders}

Given a graph $G$, let $\spem{L}[G]$ denote the vector space with basis given by all linear orders on $V(G)$.
We give $\spem{L}$ the structure of a $(q,t)$-Hopf monoid, denoted $\spem{L}_{q,t}$.
Given disjoint sets $S$ and $T$, with linear orders $\ell_S$ and $\ell_T$, let $\ell_S \cdot \ell_T$ denote the concatenation, which is a linear order on $S \sqcup T$ such that $i < j$ in the order if either $i \in S, j \in T$, or $\{i,j \} \subset R \in \{S, T \}$, and $i < j$ in $\ell_R$. Also, given a linear order $\ell$ on $I$, let $\ell|_I$ denote the restriction of the linear order to $I$. Let $inv_{S,T}(\ell, G) = |\{st \in E(G): s \in S, t \in T, s \geq_{\ell} t \}|$. Finally, let $\bar{\ell}$ denote the dual linear order.

\begin{proposition}
 $\spem{L}_{q,t}$ is a $(q,t)$-Hopf monoid. The $(G,S,T)$-component of the product is given by:
$$\spem{L}_{q,t}[G_S] \otimes \spem{L}_{q,t}[G_T] \to \spem{L}_{q,t}[G]$$
$$\ell_S \otimes \ell_T \mapsto \ell_S \cdot \ell_T$$
The $(G,S,T)$-component of the coproduct is given by:
$$\spem{L}_{q,t}[G] \to \spem{L}_{q,t}[G_S] \otimes \spem{L}_{q,t}[G_T]$$
$$\ell \mapsto q^{inv_{S,T}(\ell, G)} t^{inv_{S,T}, \overline{G}} \ell_S \otimes \ell_T$$
The $G$-component of the antipode is given by:
$$s_G: \ell \mapsto (-1)^n q^{e(G)} t^{e(\overline{G})} \bar{\ell}$$
\end{proposition}

\subsection{Acyclic orientations}

Let $\overrightarrow{\spem{A}}[G]$ denote the vector space with basis $\overrightarrow{O}$, indexed by acyclic orientations of $G$.
It is not hard to see that $\overrightarrow{\spem{A}}$ is a graphical species.
Acyclic orientations arise in the study of graphic arrangements \cite{acyclic-graphs2}, and the study of graph coloring \cite{acyclic-graphs}. Let $\chi(G, q)$ denote the chromatic polynomial of $G$. Then Stanley \cite{acyclic-graphs} showed that $(-1)^{|G|}\chi(G,-1)$ is the number of acyclic orientations of $G$. There is also a bijection between the acyclic orientations of $G$ and the regions of the graphic hyperplane arrangement \cite{acyclic-graphs2}. In this sense, acyclic orientations arise as a geometrically-motivated graph-analogue of linear orders.

We proceed to give $\overrightarrow{\spem{A}}$ the structure of a $q$-Hopf monoid. Let $G$ be a finite graph, and let $S|T \models V(G)$. 
Let $\overrightarrow{O}$ be an acyclic orientation of $G(S)$, and $\overrightarrow{P}$ be an acyclic orientation of $G(T)$. 
Then we obtain $\overrightarrow{O} \cdot \overrightarrow{P}$ from $\overrightarrow{O} \cup \overrightarrow{P}$ by adding directed edges $(s,t)$, whenever $s \in S, t \in T$ and $st \in E(G)$. This defines an acyclic orientation on $G$.
Also, $e_{S,T}(\overrightarrow{O})$ is the number of directed edges in $\overrightarrow{O}$ of the form $(s,t)$ where $s \in S, t \in T$.

Given an acyclic orientation $\overrightarrow{O}$ on $G$, and a subset $S \subset V(G)$, $\overrightarrow{O}|_S$ is acyclic orientation with vertex set $S$, and directed edges $(u,v)$ whenever $u,v \in S$ and $(u,v) \in E(\overrightarrow{O})$.
Finally, $\overleftarrow{O}$ is the acyclic orientation of $G$ with edges $(u,v)$ if and only if $(v,u)$ is an edge of $\overrightarrow{O}$.

\begin{proposition}
 $\overrightarrow{\spem{A}}_q$ is a cocommutative, connected $q$-Hopf monoid. The $(G, S, T)$-component of the product is given by
$$\mu_G^{S,T}: \overrightarrow{\spem{A}}_{q}[G(S)] \otimes \overrightarrow{\spem{A}}_q[G(T)] \to \overrightarrow{\spem{A}}_q[G]$$
$$ \overrightarrow{O} \otimes \overrightarrow{P} \mapsto \overrightarrow{O} \cdot \overrightarrow{P}$$
The $(G,S,T)$-component of the coproduct is given by:
$$\Delta_G^{S,T}: \overrightarrow{\spem{A}}_q[G] \to \overrightarrow{\spem{A}}_q[G(S)] \otimes \overrightarrow{\spem{A}}_q[G(T)]$$
$$\overrightarrow{O} \mapsto q^{e_{T, S}(\overrightarrow{O})} \overrightarrow{O}|_S \otimes \overrightarrow{O}|_T$$
The $G$-component of the antipode is given by:
$$s_G: \overrightarrow{O} \mapsto (-1)^{|V(G)|} q^{e_G} \overleftarrow{O}$$

Moreover, $\overrightarrow{\spem{A}}$ is disjoint commutative.
\end{proposition}

Let $G$ be the following graph:
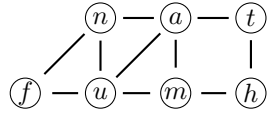
\begin{center}
\begin{tikzpicture}[>=latex]
 \node[circle] (f) at (0,-1) {$f$};
\draw (0,-1) circle (0.2cm);
\node[circle] (u) at (1,-1) {$u$};
\draw (1,-1) circle (0.2cm);
\node[circle] (n) at (1,0) {$n$};
\draw (1,0) circle (0.2cm);
\node[circle] (m) at (2,-1) {$m$};
\draw (2,-1) circle (0.2cm);
\node[circle] (a) at (2,0) {$a$};
\draw (2,0) circle (0.2cm);
\node[circle] (t) at (3,0) {$t$};
\draw (3,-1) circle (0.2cm);
\node[circle] (h) at (3,-1) {$h$};
\draw (3,0) circle (0.2cm);

\draw[-, thick] (f) -- (u);
\draw[-, thick] (f) -- (n);
\draw[-, thick] (u) -- (n);
\draw[-, thick] (m) -- (a);
\draw[-, thick] (a) -- (t);
\draw[-, thick] (t) -- (h);
\draw[-, thick] (m) -- (h);
\draw[-, thick] (a) -- (n);
\draw[-, thick] (m) -- (u);
\draw[-, thick] (a) -- (u);

\end{tikzpicture}
\end{center}
An example of the $(G, \{f,u,n\}, \{m,a,t,h\})$-component of the product is given by:
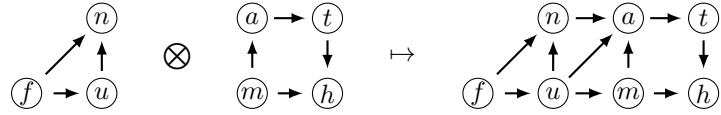
\begin{center}
\begin{tikzpicture}[>=latex]
 \node[circle] (f) at (0,-1) {$f$};
\draw (0,-1) circle (0.2cm);
\node[circle] (u) at (1,-1) {$u$};
\draw (1,-1) circle (0.2cm);
\node[circle] (n) at (1,0) {$n$};
\draw (1,0) circle (0.2cm);
\node[circle] (m) at (3,-1) {$m$};
\draw  (3,-1) circle (0.2cm);
\node[circle] (a) at (3,0) {$a$};
\draw  (3,0) circle (0.2cm);
\node[circle] (t) at (4,0) {$t$};
\draw  (4,-1) circle (0.2cm);
\node[circle] (h) at (4,-1) {$h$};
\draw  (4,0) circle (0.2cm);
 \node[circle] (f1) at (6,-1) {$f$};
\draw (6,-1) circle (0.2cm);
\node[circle] (u1) at (7,-1) {$u$};
\draw (7,-1) circle (0.2cm);
\node[circle] (n1) at (7,0) {$n$};
\draw (7,0) circle (0.2cm);
\node[circle] (m1) at (8,-1) {$m$};
\draw  (8,-1) circle (0.2cm);
\node[circle] (a1) at (8,0) {$a$};
\draw  (8,0) circle (0.2cm);
\node[circle] (t1) at (9,0) {$t$};
\draw  (9,-1) circle (0.2cm);
\node[circle] (h1) at (9,-1) {$h$};
\draw  (9,0) circle (0.2cm);
\node[circle] (ot) at (2, -0.5) {$\bigotimes$};
\node[circle] (ra) at (5, -0.5) {$\mapsto$};

\draw[->, thick] (f) -- (u);
\draw[->, thick] (f) -- (n);
\draw[->, thick] (u) -- (n);
\draw[->, thick] (m) -- (a);
\draw[->, thick] (a) -- (t);
\draw[->, thick] (t) -- (h);
\draw[->, thick] (m) -- (h);
\draw[->, thick] (f1) -- (u1);
\draw[->, thick] (f1) -- (n1);
\draw[->, thick] (u1) -- (n1);
\draw[->, thick] (m1) -- (a1);
\draw[->, thick] (a1) -- (t1);
\draw[->, thick] (t1) -- (h1);
\draw[->, thick] (m1) -- (h1);
\draw[->, thick] (n1) -- (a1);
\draw[->, thick] (u1) -- (m1);
\draw[->, thick] (u1) -- (a1);
\end{tikzpicture}
\end{center}

An example of the $(G, \{f,u,n\}, \{m,a,t,h\})$-component of the coproduct is given by:
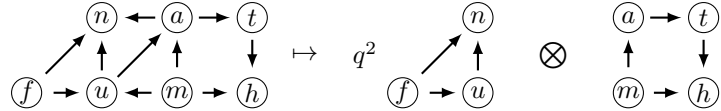
\begin{center}
\begin{tikzpicture}[>=latex]
 \node[circle] (f) at (0,-1) {$f$};
\draw (0,-1) circle (0.2cm);
\node[circle] (u) at (1,-1) {$u$};
\draw (1,-1) circle (0.2cm);
\node[circle] (n) at (1,0) {$n$};
\draw (1,0) circle (0.2cm);
\node[circle] (m) at (2,-1) {$m$};
\draw  (2,-1) circle (0.2cm);
\node[circle] (a) at (2,0) {$a$};
\draw  (2,0) circle (0.2cm);
\node[circle] (t) at (3,0) {$t$};
\draw  (3,-1) circle (0.2cm);
\node[circle] (h) at (3,-1) {$h$};
\draw  (3,0) circle (0.2cm);
 \node[circle] (f1) at (5,-1) {$f$};
\draw (5,-1) circle (0.2cm);
\node[circle] (u1) at (6,-1) {$u$};
\draw (6,-1) circle (0.2cm);
\node[circle] (n1) at (6,0) {$n$};
\draw (6,0) circle (0.2cm);
\node[circle] (m1) at (8,-1) {$m$};
\draw  (8,-1) circle (0.2cm);
\node[circle] (a1) at (8,0) {$a$};
\draw  (8,0) circle (0.2cm);
\node[circle] (t1) at (9,0) {$t$};
\draw  (9,-1) circle (0.2cm);
\node[circle] (h1) at (9,-1) {$h$};
\draw  (9,0) circle (0.2cm);
\node[circle] (ot) at (7, -0.5) {$\bigotimes$};
\node[circle] (ra) at (3.7, -0.5) {$\mapsto$};
\node[circle] (q) at (4.5, -0.5) {$q^2$};

\draw[->, thick] (f) -- (u);
\draw[->, thick] (f) -- (n);
\draw[->, thick] (u) -- (n);
\draw[->, thick] (m) -- (a);
\draw[->, thick] (a) -- (t);
\draw[->, thick] (t) -- (h);
\draw[->, thick] (m) -- (h);
\draw[->, thick] (a) -- (n);
\draw[->, thick] (m) -- (u);
\draw[->, thick] (u) -- (a);
\draw[->, thick] (f1) -- (u1);
\draw[->, thick] (f1) -- (n1);
\draw[->, thick] (u1) -- (n1);
\draw[->, thick] (m1) -- (a1);
\draw[->, thick] (a1) -- (t1);
\draw[->, thick] (t1) -- (h1);
\draw[->, thick] (m1) -- (h1);
\end{tikzpicture}
\end{center}

Finally, an example of the $G$-component of the antipode is given by:
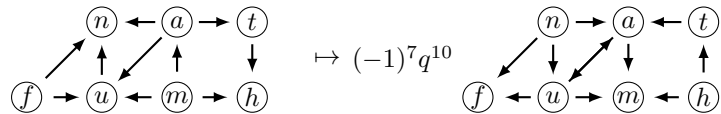
\begin{center}
\begin{tikzpicture}[>=latex]
 \node[circle] (f) at (0,-1) {$f$};
\draw (0,-1) circle (0.2cm);
\node[circle] (u) at (1,-1) {$u$};
\draw (1,-1) circle (0.2cm);
\node[circle] (n) at (1,0) {$n$};
\draw (1,0) circle (0.2cm);
\node[circle] (m) at (2,-1) {$m$};
\draw (2,-1) circle (0.2cm);
\node[circle] (a) at (2,0) {$a$};
\draw (2,0) circle (0.2cm);
\node[circle] (t) at (3,0) {$t$};
\draw (3,-1) circle (0.2cm);
\node[circle] (h) at (3,-1) {$h$};
\draw (3,0) circle (0.2cm);

\draw[->, thick] (f) -- (u);
\draw[->, thick] (f) -- (n);
\draw[->, thick] (u) -- (n);
\draw[->, thick] (m) -- (a);
\draw[->, thick] (a) -- (t);
\draw[->, thick] (t) -- (h);
\draw[->, thick] (m) -- (h);
\draw[->, thick] (a) -- (n);
\draw[->, thick] (m) -- (u);
\draw[->, thick] (a) -- (u);
\node[circle] (ma) at (4, -0.5) {$\mapsto$};
\node[circle] (q) at (5, -0.5) {$(-1)^7 q^{10}$};
 \node[circle] (f1) at (6,-1) {$f$};
\draw (6,-1) circle (0.2cm);
\node[circle] (u1) at (7,-1) {$u$};
\draw (7,-1) circle (0.2cm);
\node[circle] (n1) at (7,0) {$n$};
\draw (7,0) circle (0.2cm);
\node[circle] (m1) at (8,-1) {$m$};
\draw (8,-1) circle (0.2cm);
\node[circle] (a1) at (8,0) {$a$};
\draw (8,0) circle (0.2cm);
\node[circle] (t1) at (9,0) {$t$};
\draw (9,-1) circle (0.2cm);
\node[circle] (h1) at (9,-1) {$h$};
\draw (9,0) circle (0.2cm);

\draw[->, thick] (u1) -- (f1);
\draw[->, thick] (n1) -- (f1);
\draw[->, thick] (n1) -- (u1);
\draw[->, thick] (a1) -- (m1);
\draw[->, thick] (t1) -- (a1);
\draw[->, thick] (h1) -- (t1);
\draw[->, thick] (h1) -- (m1);
\draw[->, thick] (n1) -- (a1);
\draw[->, thick] (u1) -- (m1);
\draw[->, thick] (u1) -- (a1);
\draw[->, thick] (a1) -- (u1);
\end{tikzpicture}
\end{center}

\begin{proof}
 One can verify that $\overrightarrow{\spem{A}}_q$ is a connected $q$-bimonoid. We prove the formula for the antipode using the Milnor-Moore formula, and induction on $|G|$. 
Let $G$ be a finite graph, $\overrightarrow{O}$ be an acyclic orientation of $G$. Recall that a sink in a directed graph is a vertex $u$, such that there is no directed edge of the form $(u,v)$. 
Clearly $\overrightarrow{O}$ must have a sink $x$. Let $S = \{x \}$, $T = V(G) \setminus \{x\}$. 
Then by induction, we see that $-\mu^G_{S,T} \circ (id_{G_S} \otimes s_{G_T}) \circ \Delta^G_{S,T}(\overrightarrow{O}) = (-1)^{|G|} q^{e_G} \overleftarrow{O}$. So it is enough to show that the other terms in the Milnor-Moore formula cancel. However, for any finite set $S$, with $x \not\in S$, a straightforward computation shows that $-\mu^G_{S,T} \circ (id_{G_S} \otimes s_{G_T}) \circ \Delta^G_{S,T}(\overrightarrow{O}) = \mu^G_{S',T'} \circ (id_{G_{S'}} \otimes s_{G_{T'}}) \circ \Delta^G_{S',T'}(\overrightarrow{O})$, where $S' = S \cup \{x\}$, $T' = T \setminus \{x\}$. Hence the other terms cancel, and by induction, we have the antipode formula.
\end{proof}

\subsection{Set compositions and stable compositions}

We define the graphical \\ species of set compositions, $\spem{\Sigma}_{q,t}$, and the graphical species of stable set compositions, $\spem{S\Sigma}_{q,t}$. Given a graph $G$, $\spem{\Sigma}[G]$ is the vector space spanned by all set compositions of $V(G)$. Stable compositions are set compositions $C = S_1|S_2|\cdots|S_k$ where $e(G_{S_i}) = 0$ for all $i$. The unordered version, stable partitions, arises in the study of chromatic symmetric functions, which we discuss later in this paper. Let $\spem{S\Sigma}[G]$ denote the vector space spanned by all stable compositions of $G$. Then $\spem{\Sigma}$ and $\spem{S\Sigma}$ are both graphical species.

Given $V = S \sqcup T$, $C = S_1 |\cdots|S_k \models S$, $C' = T_1|\cdots|T_m \models T$, define $C \cdot C' = S_1|\cdots|S_k|T_1|\cdots|T_m$. Given a set composition $C = V_1|\cdots|V_k \models V(G)$, define $C_S$ to be obtained from $V_1 \cap S|\cdots|V_k \cap S$ by removing all empty blocks. 

Given $C = C_1|\cdots|C_k \models V(G)$, $S|T \models V(G)$, let $inv_{S,T}(C,G) = |\{(s,t): s \in S, t \in T, st \in E(G), s \in C_j, t \in C_i, j > i \}|$. Finally, $e(G,C) = |\{uv: uv \in E(G), u \in C_i, v \in C_j, i \neq j \}|$, and $|C| = k$.

\begin{theorem}
$\spem{\Sigma}_{q,t}$ is a $(q,t)$-Hopf monoid. The $(G,S,T)$-component of the product is given by $$\mu_G^{S,T}: \spem{\Sigma}[G_S] \otimes \spem{\Sigma}[G_T] \to \spem{\Sigma}[G]$$
$$C \otimes C' \mapsto C \cdot C'$$
The $(G,S,T)$-component of the coproduct is given by $$\Delta_G^{S,T}: \spem{\Sigma}[G] \to \spem{\Sigma}[G_S] \otimes \spem{\Sigma}[G_T]$$ $$C \mapsto q^{inv_{S,T}(C,G,)} t^{inv_{S,T}(C, \overline{G})} C_S \otimes C_T$$
The $G$-component of the antipode is given by $$S_G: \spem{\Sigma}[G] \to \spem{\Sigma}[G]$$ $$C \mapsto \sum_{C' \leq \overline{C}} (-1)^{|C'|} q^{e(G,C')} t^{e(\overline{G}, T)} C'$$

\end{theorem}

 It turns out that $\spem{S\Sigma}_{q,t}$ is a subHopf monoid of $\spem{\Sigma}$. That is, the product, coproduct, and antipode for $\spem{\Sigma}_{q,t}$ also define the Hopf monoid structure for $\spem{S\Sigma}_{q,t}$.

\begin{theorem}
$\spem{S\Sigma}_{q,t}$ is a subHopf monoid of $\spem{\Sigma}_{q,t}$
\end{theorem}

\subsection{Relating the various graphical species}

There is already a known Hopf monoid morphism between $\spem{L}$ and $\spem{\Sigma}$ in the category of species. The map is given by sending the linear order $\ell$ to the set composition $C(\ell) = \ell_1|\ell_2|\cdots|\ell_k$. That is, we view each element as a singleton block, and linear order the blocks according to $\ell$. Let us denote the resulting morphism by $\iota$. There is also a well-known morphism from $\spem{L}$ to $\spem{E}$. We show  how these two morphisms factor through $\spem{S\Sigma}$ and $\overrightarrow{\spem{A}}$.

Clearly $\iota$ is a morphism of graphical species. Note that $C(\ell)$ is a stable composition, so we also obtain a map $\iota: \spem{L}_{q,t} \to \spem{S\Sigma}_{q,t}$.

\begin{proposition}
The following is a commutative diagram of Hopf monoids and Hopf monoid morphisms. The maps are all inclusions.
\begin{center}
\begin{tikzpicture}[>=latex]
\node (a) at (0,0) {$\spem{L}_{q,t}$};
 \node (b) at (0,-2) {$\spem{S\Sigma}_{q,t}$};
 \node (c) at (2,-2) {$\spem{\Sigma}_{q,t}$};

 \draw[->] (a) --  node[left] {$\iota$} (b);
 \draw[->] (b) -- node[below] {$\iota$} (c);
 \draw[->] (a) -- node[above] {$\iota$} (c);
\end{tikzpicture}
\end{center}

\end{proposition}

Next we show that the morphism $\iota: \spem{L}_{q,t} \to \spem{S\Sigma}_{q,t}$ is part of a commutative diagram involving $\overrightarrow{\spem{A}}_q$ when $t = 1$.

We construct maps $\overrightarrow{\pi}$, from $\spem{L}_q$, and $\spem{S\Sigma}_q$, to $\overrightarrow{\spem{A}}_q$.
Given a graph $G$, and a linear order $\ell$, there is a natural acyclic orientation $\overrightarrow{O}(\ell)$ associated to $\ell$: we direct the edge $uv$ from $u$ to $v$ if $u <_{\ell} v$. This defines a map $\overrightarrow{\pi}: \spem{L}_{q,1} \to \overrightarrow{\spem{A}}$, which is also a Hopf monoid morphism. 
Finally, given a stable composition $C = S_1|\cdots|S_k$, there is a natural acyclic orientation $\overrightarrow{O}(C)$. One directs an edge $uv$ from $u$ to $v$ if $u \in S_i$, $v \in S_j$ and $i < j$. Note that by definition of stable composition, $u$ and $v$ must be in different blocks. This construction gives us a third map $\overrightarrow{\pi}: \spem{S}\Sigma_{q,1} \to \overrightarrow{\spem{A}}_q$.

\begin{proposition}
We have the following commutative diagram of Hopf monoids and Hopf monoid morphisms:
\begin{center}
\begin{tikzpicture}[>=latex]
\node (a) at (0,0) {$\spem{L}_q$};
 \node (b) at (0,-2) {$\spem{S\Sigma}_{q}$};
 \node (c) at (2,-2) {$\overrightarrow{\spem{A}}_{q}$};

 \draw[->] (a) --  node[left] {$\iota$} (b);
 \draw[->] (b) -- node[below] {$\overrightarrow{\pi}$} (c);
 \draw[->] (a) -- node[above] {$\overrightarrow{\pi}$} (c);
\end{tikzpicture}
\end{center}
\end{proposition}

Finally, recall the map $\pi: \spem{L} \to \spem{E}$, whose $G$-component is given by:
$$\pi_G: \spem{L}[G] \to \spem{E}[G]$$
$$ \ell \mapsto 1$$
This is known as the abelianization map. There is a similar map $\pi: \overrightarrow{\spem{A}} \to \spem{E}$, whose $G$-component is also the surjection which maps every acyclic orientation to $1$.

The final result of this section is that the abelianization map factors.
\begin{proposition}
We have the following commutative diagram of Hopf monoids and Hopf monoid morphisms:
\begin{center}
\begin{tikzpicture}[>=latex]
\node (a) at (0,0) {$\spem{L}$};
 \node (b) at (0,-2) {$\overrightarrow{\spem{A}}$};
 \node (c) at (2,-2) {$\spem{E}$};

 \draw[->] (a) --  node[left] {$\overrightarrow{\pi}$} (b);
 \draw[->] (b) -- node[below] {$\pi$} (c);
 \draw[->] (a) -- node[above] {$\pi$} (c);
\end{tikzpicture}
\end{center}
\end{proposition}

As we shall see in the next section, it is not the only way to factor the abelianization map. 

\section{Graphical Species generalizing $\spem{\Pi}$}

In this section, we construct three graphical species which generalize the species of set partitions. The first, of course, is the graphical species of vertex set partitions $\spem{\Pi}$. The second is $\spem{S\Pi}$, the graphical species of stable partitions. Finally, we have $\spem{FL}$, the graphical species of flats.

\subsection{Partitions and stable set partitions}

For a graph $G$, let $\spem{\Pi}[G]$ denote the vector space with basis $\{m_{\pi}: \pi \vdash V(G) \}$. It is clear that $\spem{\Pi}$ is a graphical species. We turn it into a Hopf monoid.

The $(G,S,T)$-component of the product is given by 
$$\mu^G_{S,T}: \spem{\Pi}[G_S] \otimes \spem{\Pi}[G_T] \to \spem{\Pi}[G]$$
$$ m_{\sigma} \otimes m_{\tau} \mapsto m_{\sigma \cup \tau}$$
The $(G,S,T)$-component of the coproduct is given by 
$$\Delta^G_{S,T}: \spem{\Pi}[G] \to \spem{\Pi}[G_S] \otimes \spem{\Pi}[G_T]$$
$$m_{\pi} \mapsto m_{\pi|_S} \otimes m_{\pi|_T}$$
where $\pi|_S$ is the restriction of $\pi$ to $S$.

Recall that set partitions are ordered by refinement.
$\spem{\Pi}$ has a second basis, the $p$ basis, given by $$m_{\pi} = \sum_{\tau \leq \pi} p_{\tau}.$$

\begin{proposition}
$\spem{\Pi}$ is a commutative, cocommutative self-dual Hopf monoid. The product in the $p$ basis 
is given by:
$$p_{\sigma} \otimes p_{\tau} \mapsto p_{\sigma \cup \tau}$$
The coproduct is given by:
\begin{displaymath}
p_{\pi} \mapsto \left\{ \begin{array}{cc} p_{\pi|_S} \otimes p_{\pi|_T} & \mbox{ if } \pi = \pi|_S \cup \pi|_T \\ 0 & \end{array} \right.
\end{displaymath}
The antipode is given by
$$s(p_{\pi}) = (-1)^{|\pi|} p_{\pi}$$
\end{proposition}

Now we construct $\spem{S\Pi}$, the graphical species of stable partitions. A stable partition of a graph $G$ is a set partition $\pi$ of $V(G)$ such that, for every edge $uv \in E(G)$, $u$ and $v$ lie in different parts of $\pi$. Let $\spem{S\Pi}[G]$ be the subspace of $\spem{\Pi[G]}$ generated by all $m_{\pi}$ where $\pi$ is a stable partition. Note that stable partitions are closed under refinement. That is, they form an order ideal in the partition lattice. Hence $\spem{\Pi[G]}$ is also the subspace generated by all $p_{\pi}$, where $\pi$ is a stable partition.

\begin{proposition}
$\spem{S\Pi}$ is a commutative, cocommutative self-dual Hopf submonoid of $\spem{\Pi}$.
\end{proposition}

Stanley introduced the chromatic symmetric function \cite{stanley-chromatic}. He expressed the chromatic symmetric function of a graph in terms of several well-known bases for symmetric functions. For the `augmented' monomial symmetric functions, the expression he found involves summing over all stable partitions. 

\subsection{Hopf monoid of flats}

Finally, we study $\spem{FL}$, the graphical species of flats. Let $G$ be a graph, and let $D \subset E(G)$. Then $(V(G), D)$ is a graph, and we can partition $V(G)$ into a set partition $\pi_D$, where two vertices belong to the same part iff they belong to the same connected component of $(V(G), D)$. Such a subset $D$ is called a \emph{flat} if $$D = \bigcup_{B \in \pi_D}  E(G_B).$$
Let $\spem{FL}[G]$ be the vector space with basis $\{M_F: F \mbox{ is a flat of } G \}$.
Then $\spem{FL}[G]$ is a graphical species. Given a subset $S$ of vertices, and a flat $F$, let $F|_S = \{uv \in F: u,v \in S \}$. Clearly, if $F$ is a flat of $G$, $F|_S$ is a flat of $G_S$. 
Flats arise naturally in the study of hyperplane arrangements and matroids. See \cite{orlik-terao} for more about hyperplane arrangements, including the intersection lattice (whose elements are flats). For more about matroids, see \cite{matroids}.

Also recall that the flats of a graph are ordered by inclusion (the so called bond lattice, or lattice of contractions of $G$). We construct another basis, the $p$ basis, by $$M_F = \sum_{F' \leq F} P_{F'}.$$

We turn $\spem{FL}[G]$ into a Hopf monoid.
\begin{proposition}
$\spem{FL}[G]$ is a commutative, cocommutative self-dual Hopf monoid. 
The $(G,S,T)$-component of the product is given by:
$$\mu_G^{S,T}: \spem{FL}[G_S] \otimes \spem{FL}[G_T] \to \spem{FL}[G]$$
$$M_F \otimes M_H \mapsto M_{F \sqcup H}$$
$$P_F \otimes P_H \mapsto P_{F \sqcup H}$$
The $(G,S,T)$-component of the coproduct is given by:
$$\Delta_G^{S,T}: \spem{FL}[G] \to \spem{FL}[G_S] \otimes \spem{FL}[G_T]$$
$$M_F \mapsto M_{F|_S} \otimes M_{F|_T}$$
\begin{displaymath}
P_F \mapsto \left\{ \begin{array}{cc} P_{F|_S} \otimes P_{F|_T} & F = F|_S \sqcup F|_T \\ 0 & \end{array} \right.
\end{displaymath}
The antipode is given by:
$$s_G: \spem{FL}[G] \to \spem{FL}[G]$$
$$M_F \mapsto \sum_{H \leq F} (-1)^{c(H)} o(H, F) M_H$$
$$P_F \mapsto (-1)^{c(F)} P_F$$
where $c(F)$ is the number of components of $(V, F)$, and $o(H,F)$ is the number of acyclic orientations of $(V, F) / H$, which is obtained from $(V,F)$ by contracting all the edges of $H$.
\end{proposition}

\begin{proposition}
$\spem{FL}$ is a subHopf monoid of $\spem{\Pi}$, under the map $$\iota: \spem{FL}[G] \to \spem{\Pi}[G]$$
$$P_{F} \mapsto P_{\pi_F}$$
\end{proposition}

For example, $$\iota_G: \spem{FL}[G] \to \spem{\Pi}[G]$$
$$P_{un, fu, fn, ma, at} \mapsto P_{fun/mat/h}$$ 
Pictorially, this map looks like
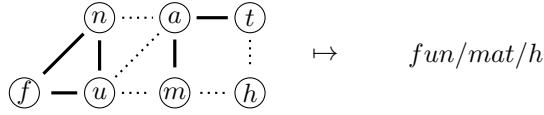
\begin{center}
\begin{tikzpicture}[>=latex]
\node (p) at (12,-0.5) {$fun/mat/h$};
\node (m) at (10, -.5) {$\mapsto$};
 \node[circle] (f1) at (6,-1) {$f$};
\draw (6,-1) circle (0.2cm);
\node[circle] (u1) at (7,-1) {$u$};
\draw (7,-1) circle (0.2cm);
\node[circle] (n1) at (7,0) {$n$};
\draw (7,0) circle (0.2cm);
\node[circle] (m1) at (8,-1) {$m$};
\draw (8,-1) circle (0.2cm);
\node[circle] (a1) at (8,0) {$a$};
\draw (8,0) circle (0.2cm);
\node[circle] (t1) at (9,0) {$t$};
\draw (9,-1) circle (0.2cm);
\node[circle] (h1) at (9,-1) {$h$};
\draw (9,0) circle (0.2cm);

\draw[-, very thick] (f1) -- (u1);
\draw[-, very thick] (f1) -- (n1);
\draw[-, very thick] (u1) -- (n1);
\draw[-, very thick] (m1) -- (a1);
\draw[-, very thick] (a1) -- (t1);
\draw[-, thick, dotted] (t1) -- (h1);
\draw[-, thick, dotted] (m1) -- (h1);
\draw[-, thick, dotted] (n1) -- (a1);
\draw[-, thick, dotted] (u1) -- (m1);
\draw[-, thick, dotted] (u1) -- (a1);
\end{tikzpicture}
\end{center}

Note that even though $\spem{FL}$ is a subHopf monoid of $\spem{\Pi}$, the inclusion map is not given by $M_F \mapsto M_{\pi_F}$. In particular, this alternate map is not even a morphism of comonoids. This is why we choose to study flats based on their edge sets, instead of viewing them as being equivalent to the notion of connected set partitions. 

Note that 
the flats of a graph $G$ correspond to the flats of the graphic arrangement for $G$. In particular, Aguiar and Mahajan studied break and join maps for flats of the braid arrangement, and if one were to generalize their work to graphic arrangements, the corresponding operations give rise to the product and coproduct of $\spem{FL}$ in the $M$ basis.

Much like $\spem{\Pi}$ has a Hopf submonoid $\spem{S\Pi}$ that has arisen in graph theory, $\spem{FL}$ also has a Hopf submonoid that has arisen in graph theory. A matching $M$ is a subset of edge of $G$ such that no two edges of $M$ have an endpoint in common. Let $\spem{M}[G]$ be the subspace of $\spem{FL}[G]$ generated by $\{M_F: F \mbox{ is a matching } \}$. Note that matchings form an order ideal in the bond lattice, hence $\spem{M}[G]$ is also generated by $\{P_F: F \mbox{ is a matching } \}$.

\begin{proposition}
$\spem{M}$ is a self-dual, commutative, cocommutative sub Hopf monoid of $\spem{FL}$.
\end{proposition}

\subsection{Commutative diagrams involving generalizations of $\spem{\Pi}$}

Now we shall relate the various generalizations of $\spem{\Pi}$ considered in this section. There is a connection between $\spem{\Pi}$ and $\spem{FL}$ involving the $M$ basis. That is, given a set partition $\pi$, let $$F(\pi) = \bigcup_{B \in \pi} E(G_B).$$
Then the map $\varphi$, whose $G$ component is
$$\varphi_G: \spem{\Pi}[G] \to \spem{FL}[G] $$
$$M_{\pi} \mapsto M_{F(\pi)}$$ 
is a Hopf monoid morphism.

For example, $$\varphi_G: \spem{\Pi}[G] \to \spem{FL}[G]$$
$$M_{un/fmat/h} \mapsto M_{un, ma, at}$$ 
Pictorially, this map looks like
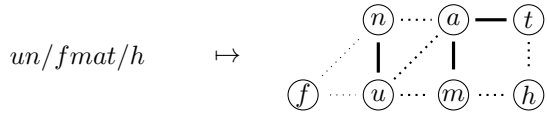
\begin{center}
\begin{tikzpicture}[>=latex]
\node (p) at (3,-0.5) {$un/fmat/h$};
\node (m) at (5, -.5) {$\mapsto$};
 \node[circle] (f1) at (6,-1) {$f$};
\draw (6,-1) circle (0.2cm);
\node[circle] (u1) at (7,-1) {$u$};
\draw (7,-1) circle (0.2cm);
\node[circle] (n1) at (7,0) {$n$};
\draw (7,0) circle (0.2cm);
\node[circle] (m1) at (8,-1) {$m$};
\draw (8,-1) circle (0.2cm);
\node[circle] (a1) at (8,0) {$a$};
\draw (8,0) circle (0.2cm);
\node[circle] (t1) at (9,0) {$t$};
\draw (9,-1) circle (0.2cm);
\node[circle] (h1) at (9,-1) {$h$};
\draw (9,0) circle (0.2cm);

\draw[-, dotted] (f1) -- (u1);
\draw[-, dotted] (f1) -- (n1);
\draw[-, very thick] (u1) -- (n1);
\draw[-, very thick] (m1) -- (a1);
\draw[-, very thick] (a1) -- (t1);
\draw[-, thick, dotted] (t1) -- (h1);
\draw[-, thick, dotted] (m1) -- (h1);
\draw[-, thick, dotted] (n1) -- (a1);
\draw[-, thick, dotted] (u1) -- (m1);
\draw[-, thick, dotted] (u1) -- (a1);
\end{tikzpicture}
\end{center}

\begin{theorem}
We have the following diagram of Hopf monoids and Hopf monoid morphisms:

\begin{center}
\begin{tikzpicture}{>=latex}
\node (a) at (0,0) {$\spem{E}$};
\node (b) at (0,2) {$\spem{S\Pi}$};
\node (d) at (2,0) {$\spem{FL}$};
\node (c) at (2,2) {$\spem{\Pi}$};
\draw[->] (b) -- node[left] {$\pi$} (a);
\draw[->] (b) -- node[above] {$\iota$} (c);
\draw[->] (a) --  node[below] {$\iota$} (d);
\draw[->] (c) -- node[right] {$\varphi$} (d);
\end{tikzpicture}
\end{center}
\end{theorem}

Let $\pi: \spem{\Sigma} \to \spem{\Pi}$ be the map which sends a composition $C$ to $\pi(C)$. 
Similarly, define $\pi: \spem{S\Sigma} \to \spem{S\Pi}$.

\begin{proposition}
We have the following commutative diagram of Hopf monoids and Hopf monoid morphisms:
\begin{center}
\begin{tikzpicture}{>=latex}
\node (a) at (0,0) {$\spem{S\Pi}$};
\node (b) at (0,2) {$\spem{S\Sigma}$};
\node (d) at (2,0) {$\spem{\Pi}$};
\node (c) at (2,2) {$\spem{\Sigma}$};
\draw[->] (b) -- node[left] {$\pi$} (a);
\draw[->] (b) -- node[above] {$\iota$} (c);
\draw[->] (a) --  node[below] {$\iota$} (d);
\draw[->] (c) -- node[right] {$\pi$} (d);
\end{tikzpicture}
\end{center}
\end{proposition}
\begin{theorem}
We have the following commutative diagram of Hopf monoids:
\begin{center}
\begin{tikzpicture}{>=latex}
\node (a) at (0,0) {$\spem{S\Pi}$};
\node (b) at (0,2) {$\spem{S\Sigma}$};
\node (d) at (2,0) {$\spem{E}$};
\node (c) at (2,2) {$\overrightarrow{\spem{A}}$};
\node (e) at (0,4) {$\spem{L}$};
\draw[->] (e) -- node[left] {$\iota$} (b);
\draw[->] (e) -- node[above] {$\overrightarrow{\pi}$} (c);
\draw[->] (b) -- node[left] {$\pi$} (a);
\draw[->] (b) -- node[above] {$\overrightarrow{\pi}$} (c);
\draw[->] (a) --  node[below] {$\rho$} (d);
\draw[->] (c) -- node[right] {$\pi$} (d);
\end{tikzpicture}
\end{center}
In particular, the resulting map from $\spem{L}$ to $\spem{E}$ is the abelianization map $\pi$.
\end{theorem}

\section{Functors on Graphical Species}

\subsection{Complementation}
We note that duality $\spem{p}^{\ast}[G] = \spem{p}[G]^{\ast}$ can be used to construct a bistrong contravariant functor. $((-)^{\ast}, \phi, \phi_0, \psi, \psi_0): (\gs, \cdot, \beta_{q,t}) \to (\gs, \cdot, \beta_{q,t})$, similar to the duality functor on species, graded vector spaces, and vector spaces. 

However, for graphical species there is a second notion of duality, known as graph complementation.
Given a finite graph $G$, $\overline{\spem{p}}[G] = \spem{p}[\overline{G}]$. This defines a functor $\overline{(-)}: \gs \to \gs$. We turn it into a bistrong monoidal functor. 

For each pair of graphical species, $\spem{g}$ and $\spem{h}$, we will define maps:

\begin{center}
\begin{tikzpicture}
    \node (A) at (0,0) {$\overline{\spem{e}} \cdot \overline{\spem{g}}$};
    \node (B) at (3,0) {$\overline{\spem{g}\cdot \spem{h}}$};
    \draw[transform canvas={yshift=0.5ex},->] (A) -- node[above] {$\varphi_{\spem{g}, \spem{h}}$} (B);
    \draw[transform canvas={yshift=-0.5ex},<-] (A) -- node[below] {$\psi_{\spem{g}, \spem{h}}$} (B);
\end{tikzpicture} 
\end{center}

Given a finite graph $G$, and $S|T \models V(G)$, we define 
\begin{center}
\begin{tikzpicture}
    \node (A) at (0,0) {$\overline{\spem{g}}[G_S] \otimes \overline{\spem{h}}[G_T] = \spem{g}[\overline{G_S}] \otimes \spem{h}[\overline{G_T}]$};
    \node (B) at (6,0) {$\spem{g}[\overline{G}_S] \otimes \spem{h}[\overline{G}_T]$};
    \draw[transform canvas={yshift=0.5ex},->] (A) -- node[above] {$\varphi_{\spem{g}, \spem{h}, G}^{S,T}$} (B);
    \draw[transform canvas={yshift=-0.5ex},<-] (A) -- node[below] {$\psi_{\spem{g}, \spem{h}, G}^{S,T}$} (B);
\end{tikzpicture} 
\end{center}
where $\varphi_{\spem{g}, \spem{h}}^{G,S,T}$ and $\psi_{\spem{g}, \spem{h}}^{G,S,T}$ are the obvious identity morphisms. 

Taking the direct sum over all $S|T \models V(G)$, we obtain maps:
\begin{center}
\begin{tikzpicture}
    \node (A) at (0,0) {$(\overline{\spem{g}} \cdot \overline{\spem{h}})[G]$};
    \node (B) at (5,0) {$(\spem{g} \cdot \spem{h})[\overline{G}] = (\overline{\spem{g} \cdot \spem{h}})[G]$};
    \draw[transform canvas={yshift=0.5ex},->] (A) -- node[above] {$\varphi_{\spem{g}, \spem{h}, G}$} (B);
    \draw[transform canvas={yshift=-0.5ex},<-] (A) -- node[below] {$\psi_{\spem{g}, \spem{h}, G}$} (B);
\end{tikzpicture} 
\end{center}

By definition, $\overline{\emptyset} = \emptyset$, so we let $\phi_0, \psi_0$ be identity isomorphisms.

\begin{theorem}
$(\overline{(-)}, \phi, \phi_0, \psi, \psi_0): (\gs, \cdot, \beta_{q,t}) \to (\gs, \cdot, \beta_{t,q})$ is a bistrong monoidal functor. Moreover, it is an involution.
\end{theorem}
Please note that the complementation functor interchanges the parameters on the braiding.
One corollary of this theorem is that $(1,q)$-Hopf monoids are equivalent to $(q,1)$-Hopf monoids, which is why we refer to $(q,1)$-Hopf monoids as $q$-Hopf monoids.

For examples, note that $\overline{\spem{S\Pi}}[G]$ has basis given by clique partitions (a clique partition is a set partition $\pi$ such that each block of $\pi$ induces a clique in $G$.

\subsection{Clique and Discrete functor}

First, we detail the Clique functor, and then the Discrete functor.
Given a finite set $I$, let $K_I$ denote the complete graph on $I$, and let $D_I$ denote the discrete graph on $I$. Recall that the complete graph has edges between every pair of vertices, and the discrete graph has no edge.

Given a graphical species $\spem{g}$, and a set $I$, let $\mathcal{K}(\spem{g})[I] = \spem{g}[K_I]$. We see that this defines a species $\mathcal{K}(\spem{g})$. Also, given a morphism $\alpha: \spem{g} \to \spem{h}$, we can define a morphism $\mathcal{K}(\alpha): \mathcal{K}(\spem{g}) \to \mathcal{K}(\spem{h})$, whose $I$-component is $\mathcal{K}(\alpha)_I = \alpha_{K_I}$. Thus, we have defined a functor $\mathcal{K}: \cat{GrSp} \to \cat{Sp}$.
We call this the Clique functor.

We turn $\mathcal{K}$ into a bistrong monoidal functor. It is not hard to see that $\mathcal{K}(\spem{g} \cdot \spem{h}) [I] \cong (\mathcal{K}(\spem{g}) \cdot \mathcal{K}(\spem{h}))[I]$. Hence, we let $\varphi$, $\psi$ be the corresponding isomorphisms. 

\begin{theorem}
$$(\mathcal{K}, \varphi, \psi): (\cat{GrSp}, \cdot, \beta_{q,t}) \to (\cat{Sp}, \cdot, \beta_q)$$ is a bistrong monoidal functor.
\end{theorem}

We can also define the discrete functor. One way to define it is as the composition of bistrong functors $\mathcal{K} \circ \overline{(-)}$. Explicitly, given a set $I$, we let $\mathcal{D}(\spem{g})[I] = \spem{g}[D_I]$. Moreover, the lax and colax structures come from the fact that $\mathcal{D}(\spem{g} \cdot \spem{h})[I] = (\mathcal{D}(\spem{g}) \cdot \mathcal{D}(\spem{h}))[I]$.
Since composition of bistrong functors is bistrong, we obtain:
\begin{corollary}
$$(\mathcal{D}, \varphi, \psi): (\cat{GrSp}, \cdot, \beta_{q,t}) \to (\cat{Sp}, \cdot, \beta_t)$$ is a bistrong monoidal functor. 
\end{corollary}

As an application, we note that $$\mathcal{K}(\spem{L}_{q,t}) \cong \mathcal{K}(\overrightarrow{\spem{A}}_q) \cong \mathcal{K}(\spem{S\Sigma}_{q,t}) \cong \spem{L}_q \cong \mathcal{D}(\spem{L}_{t,q})$$
$$\mathcal{K}(\spem{\Sigma}_{q,t}) \cong \spem{\Sigma}_q \cong \mathcal{D}(\spem{\Sigma}_{t,q}) \cong \mathcal{D}(\spem{S\Sigma}_{t,q})$$
$$\mathcal{K}(\spem{\Pi}) \cong \mathcal{K}(\spem{FL}) \cong \spem{\Pi} \cong \mathcal{D}(\spem{\Pi}) \cong \mathcal{D}(\spem{S\Pi})$$
$$\mathcal{K}(\spem{E}) \cong \mathcal{K}(\spem{S\Pi}) \cong \spem{E} \cong \mathcal{D}(\spem{E}) \cong \mathcal{D}(\spem{FL}) \cong \mathcal{D}(\overrightarrow{\spem{A}}_q)$$

Thus, we have shown how the examples of species given in Section 2 are all images of Hopf monoids in graphical species under the clique and discrete functors. Note that under these functors, the resulting species are all Hopf monoids. In particular, the resulting Hopf monoids are similar to $\spem{L}_{q,q}, \spem{\Sigma}_{q,q}$, $\spem{\Pi}$ and $\spem{E}$.

\section{Future Directions}

First, we have the following conjecture:
\begin{conjecture}
$\overrightarrow{\spem{A}} \circ \spem{g}$ is the free disjoint commutative Hopf monoid on a positive graphical species $\spem{g}$.
\end{conjecture}
This conjecture is actually true, but it demonstrates the fact that there are new universal constructions in the category of graphical species. There are several others we intend to present in a future paper. Most of this work will give an alternate explanation to the commutative diagrams shown in the current paper.

Of course, with a notion of substitution product, we can talk about monoids in $(\gs, \circ, X)$. We shall call such monoids `graph operads', since symmetric operads are monoids in $(\spe, \circ, X)$. Does $\overrightarrow{\spem{A}}$ have the structure of a graph operad? Does this explain the notion of disjoint commutativity? Are there other interesting graph operads? Can graph operads find application outside of graph theory? Also, the notion of graph quotient makes more sense for multigraphs. What do multigraph species look like?

There are other functors from graphical species to species. That is, given a graphical species $\spem{g}$, and a set $I$, let $$\mathcal{G}(\spem{g})[I] = \bigoplus_{G: V(G) = I} \spem{g}[G].$$ It turns out that this defines a functor $\mathcal{G}: \gs \to \spe$. This functor can be turned into a bilax monoidal functor in a variety of ways, at least one of which turns $\mathcal{G}(\spem{E})$ into the Hopf monoids $\spem{G}$ studied by Marcelo and Aguiar \cite{thebook}. These bilax monoidal functors need to be fully detailed. In particular, applying $\mathcal{G}$ to $\overrightarrow{\spem{A}}$ gives us the species of directed acyclic graphs, $\overrightarrow{\spem{D}}$ and each bilax monoidal structure on $\mathcal{G}$ gives a new Hopf monoid structure on $\overrightarrow{\spem{D}}$. These should be studied in more detail.

Finally, several Hopf algebras can be associated to sequences of polytopes: the Malvenuto-Reutenauer Hopf algebra of permutations $MR$ \cite{malvenuto} are associated to vertices of permutohedra, the Loday-Ronco Hopf algebra of planar binary trees $LR$ \cite{loday-ronco} are associated to vertices of associahedra. Moreover, the faces of these polytopes give rise to Hopf algebras as well \cite{polytopes}. Carr and Devadoss \cite{graphassoc} introduced graph associahedra, generalizing both the associahedron and the permutohedron. The vertices of such polytopes are indexed by `tubings' of $G$. In particular, Forcey and Springfield \cite{cyclohedron} have described the Hopf algebra structure of  $MR$ and $LR$ in terms of `tubings' of complete graphs, and paths, respectively. Naturally, the vertices (and faces) of graph associahedra are examples of graphical species. Do they form Hopf monoids? Does this Hopf monoid structure reflect the work of Forcey and Springfield? Can we use this information to construct new Hopf algebras coming from other sequences of graph associahedra?

\subsection{acknowledgments}

The author was partially supported by an NSF grant DMS-0932078,
administered by the Mathematical Sciences Research Institute while the
author was in residence at MSRI during the Complementary Program, Fall 2010-Spring 2011.
This work was developed during the
visits of the author to MSRI and we thank the institute for its
hospitality.

\bibliographystyle{amsalpha}

\bibliography{graphmonoids}

\providecommand{\bysame}{\leavevmode\hbox to3em{\hrulefill}\thinspace}
\providecommand{\MR}{\relax\ifhmode\unskip\space\fi MR }
\providecommand{\MRhref}[2]{%
  \href{http://www.ams.org/mathscinet-getitem?mr=#1}{#2}
}
\providecommand{\href}[2]{#2}
\begin{thebibliography}{BLL98}

\bibitem[AM10]{thebook}
Marcelo Aguiar and Swapneel Mahajan, \emph{Monoidal functors, species and
  {H}opf algebras}, CRM Monograph Series, vol.~29, American Mathematical
  Society, Providence, RI, 2010, With forewords by Kenneth Brown and Stephen
  Chase and Andr{\'e} Joyal.

\bibitem[BC10]{hyperoctahedral}
N.~Bergeron and P.~Choquette, \emph{Hyperoctahedral species}, S\'em. Lothar.
  Combin. \textbf{61A} (2009/10), Art. B61Aj, 22. \MR{2734182 (2011j:05369)}

\bibitem[B{\'e}n63]{benabou}
Jean B{\'e}nabou, \emph{Cat\'egories avec multiplication}, C. R. Acad. Sci.
  Paris \textbf{256} (1963), 1887--1890. \MR{0148719 (26 \#6225)}

\bibitem[BLL98]{species-book}
F.~Bergeron, G.~Labelle, and P.~Leroux, \emph{Combinatorial species and
  tree-like structures}, Encyclopedia of Mathematics and its Applications,
  vol.~67, Cambridge University Press, Cambridge, 1998, Translated from the
  1994 French original by Margaret Readdy, With a foreword by Gian-Carlo Rota.
  \MR{1629341 (2000a:05008)}

\bibitem[CD06]{graphassoc}
Michael~P. Carr and Satyan~L. Devadoss, \emph{Coxeter complexes and
  graph-associahedra}, Topology Appl. \textbf{153} (2006), no.~12, 2155--2168.

\bibitem[Cha00]{polytopes}
Fr{\'e}d{\'e}ric Chapoton, \emph{Alg\`ebres de {H}opf des permutah\`edres,
  associah\`edres et hypercubes}, Adv. Math. \textbf{150} (2000), no.~2,
  264--275.

\bibitem[FS10]{cyclohedron}
S.~Forcey and Derriell Springfield, \emph{Geometric combinatorial algebras:
  cyclohedron and simplex}, J. Algebraic Combin. \textbf{32} (2010), no.~4,
  597--627.

\bibitem[Gre77]{acyclic-graphs2}
Curtis Green, \emph{Acyclic orientations}, Higher Comb. (1977), 65--68.

\bibitem[Joy81]{joyal}
Andr{\'e} Joyal, \emph{Une th\'eorie combinatoire des s\'eries formelles}, Adv.
  in Math. \textbf{42} (1981), no.~1, 1--82. \MR{633783 (84d:05025)}

\bibitem[JR79]{joni-rota}
S.~A. Joni and G.-C. Rota, \emph{Coalgebras and bialgebras in combinatorics},
  Stud. Appl. Math. \textbf{61} (1979), no.~2, 93--139.

\bibitem[JS93]{joyal-street}
Andr{\'e} Joyal and Ross Street, \emph{Braided tensor categories}, Adv. Math.
  \textbf{102} (1993), no.~1, 20--78. \MR{1250465 (94m:18008)}

\bibitem[LR98]{loday-ronco}
Jean-Louis Loday and Mar{\'{\i}}a~O. Ronco, \emph{Hopf algebra of the planar
  binary trees}, Adv. Math. \textbf{139} (1998), no.~2, 293--309.

\bibitem[MR95]{malvenuto}
Clauda Malvenuto and Christophe Reutenauer, \emph{Duality between
  quasi-symmetric functions and the {S}olomon descent algebra}, J. Algebra
  \textbf{177} (1995), no.~3, 967--982.

\bibitem[OT92]{orlik-terao}
Peter Orlik and Hiroaki Terao, \emph{Arrangements of hyperplanes}, Grundlehren
  der Mathematischen Wissenschaften [Fundamental Principles of Mathematical
  Sciences], vol. 300, Springer-Verlag, Berlin, 1992.

\bibitem[Oxl92]{matroids}
James~G. Oxley, \emph{Matroid theory}, Oxford Science Publications, The
  Clarendon Press Oxford University Press, New York, 1992. \MR{1207587
  (94d:05033)}

\bibitem[Sta73]{acyclic-graphs}
Richard~P. Stanley, \emph{Acyclic orientations of graphs}, Discrete Math.
  \textbf{5} (1973), 171--178.

\bibitem[Sta95]{stanley-chromatic}
\bysame, \emph{A symmetric function generalization of the chromatic polynomial
  of a graph}, Adv. Math. \textbf{111} (1995), no.~1, 166--194. \MR{1317387
  (96b:05174)}

\end{thebibliography}

\end{document}